\documentclass{amsart}

\usepackage{amsfonts}
\usepackage{amssymb}
\newcommand{\N}{{\mathbb N}}
\newcommand{\Z}{{\mathbb Z}}

\newcommand{\ga}{$\Gamma(G,{\mathcal A})$}
\newcommand{\cA}{\mathcal{A}}
\newcommand{\cB}{\mathcal{B}}

\theoremstyle{plain}
\newtheorem{lemma}{Lemma}[section]
\newtheorem{cor}{Corollary}
\newtheorem{thm}{Theorem}

\theoremstyle{remark}
\newtheorem{remark}{Remark}[section]

\begin{document}

\title{On Residual Properties of Word Hyperbolic Groups}

\author{Ashot Minasyan}
\address{Universit\'{e} de Gen\`{e}ve,
Section de Math\'{e}matiques,
2-4 rue du Li\`{e}vre,
Case postale 64,
1211 Gen\`{e}ve 4, Switzerland}
\email{aminasyan@gmail.com}

\begin{abstract} For a fixed word hyperbolic group we compare different residual properties
related to quasiconvex subgroups.
\end{abstract}
\thanks{This work was partially supported by the NSF grant DMS \#0245600 of A. Ol'shanskii and M. Sapir.}

\maketitle

\section{Introduction}
Any group $G$ can be equipped with a {\it profinite topology}
$\mathcal{PT}(G)$, whose basic open sets are cosets to normal
finite index subgroups. It is easy to see that the group
operations are continuous in $\mathcal{PT}(G)$. The group is
residually finite if and only if the profinite topology is
Hausdorff.

A subgroup $H \le G$ is closed in $\mathcal{PT}(G)$ if and only if it is equal to an intersection of finite
index subgroups; equivalently, for any element $g \notin H$ there exists a homomorphism $\varphi$ from $G$ to a finite
group $L$ such that  $\varphi(g) \notin \varphi(H)$. In this case the subgroup $H$ is called $G$-{\it separable}.

The profinite closure of a subgroup $H \le G$, i.e., the smallest closed subset containing $H$, is equal to the intersection
of all finite index subgroups $K$ of $G$ such that $H \le K$.

A group $G$ is said to be LERF if every finitely generated subgroup is closed in $\mathcal{PT}(G)$. The class of
all LERF groups includes free groups \cite{Hall}, surface groups \cite{Scott} and fundamental groups of certain
$3$-manifolds \cite{Scott}, \cite{Gitik-LERFMan}.

Let $G$ be a (word) hyperbolic group with a finite symmetrized generating set $\cA$ and let {\ga} be the corresponding
Cayley graph of $G$. A subset $Q \subseteq G$ is said to be {\it quasiconvex} if there exists a constant
$\eta \ge 0$ such that for any pair of elements $u,v \in Q$ and any geodesic segment $p$ connecting $u$ and $v$,
$p$ belongs to a closed $\eta$-neighborhood of the subset $Q$ in {\ga}. Quasiconvex subgroups are precisely
those finitely generated subgroups which are embedded in $G$ without distortion \cite[Lemma 1.6]{hyp}.

As it was noted in \cite{prof-hyp}, in the context of word hyperbolic groups instead of studying LERF-groups
it makes sense to study GFERF-groups. A hyperbolic group is called GFERF if each quasiconvex subgroup
is closed in $\mathcal{PT}(G)$. Thus, within the class of hyperbolic groups, the notion of GFERF is more general
than LERF: every LERF hyperbolic group is GFERF but not vice versa.

Unfortunately, it is absolutely unclear how to
decide if a random word hyperbolic group is GFERF (or LERF).
With this purpose, D. Long \cite{Long-engulf} and, later, G. Niblo and B. Williams \cite{Niblo-Wil}
suggested to utilize the {\it engulfing} property. They say that a subgroup $H \le G$
is {\it engulfed} if it is contained in a proper finite index subgroup of $G$.

The author was mainly interested in the following two theorems established by Niblo and Williams in 2002:

\vspace{.3cm}
\noindent  {\bf Theorem A.} (\cite[Thm. 4.1]{Niblo-Wil}) {\it Let $G$ be a word hyperbolic group and suppose that $G$
engulfs every finitely generated
free subgroup with limit set a proper subset of the boundary of $G$. Then the intersection of all finite index
subgroups of $G$ is finite. If $G$ is torsion-free then it is residually finite.}

\vspace{.3cm}

\noindent {\bf Theorem B.} (\cite[Thm. 5.2]{Niblo-Wil}) {\it Let $G$ be a word hyperbolic group which
engulfs every finitely generated subgroup $K$ such that the limit set $\Lambda(K)$ is a proper subset
of the boundary of $G$. Then every quasiconvex subgroup of $G$ has a finite index in its profinite closure
in $G$.}

\vspace{.3cm}

The main goal of this paper is to generalize Theorems A and B by weakening their assumptions and,
in certain situations, strengthening their conclusions.

In a hyperbolic group $G$ the structure of a distorted subgroup can be very complicated. Thus, the basic idea
is to use assumptions which concern only quasiconvex subgroups. We prove the following
results in Section \ref{sec:engulf}:

\begin{thm} \label{thm:free-engulf} Let $G$ be a hyperbolic group with a generating set
of cardinality $s \in \N$. Suppose that each proper free quasiconvex subgroup of rank $s$ is engulfed
in $G$. Then the intersection of all finite index subgroups of $G$ is finite. If $G$ is torsion-free then
it is residually finite.
\end{thm}

\begin{thm} \label{thm: engulf-alm-gferf} Suppose $G$ is a hyperbolic group which engulfs each
proper quasiconvex subgroup. Let $H$ be an arbitrary quasiconvex subgroup of $G$. Then $H$ has a finite
index in its profinite closure $K$. Moreover, $K = HQ$ where $Q$ is the intersection of all finite index
subgroups of $G$.
\end{thm}

The assumptions of Theorems \ref{thm:free-engulf} and \ref{thm: engulf-alm-gferf} are, indeed, less restrictive
than the assumptions of Theorems A and B, because if $H$ is a quasiconvex subgroup of a hyperbolic group $G$ with
$|G:H|=\infty$ then the limit set $\Lambda(H)$ is a proper subset of $\partial G$ (\cite[Thm. 4]{Swenson},
\cite[Lemma 8.2]{paper2}).

In the residually finite case Theorem \ref{thm: engulf-alm-gferf} can be reformulated as follows:

\begin{thm}\label{thm:qc-englulf_2} Let $G$ be a residually finite hyperbolic group where
every proper quasiconvex subgroup is engulfed. Then $G$ is GFERF.
\end{thm}

Combining together the claims of Theorems \ref{thm:free-engulf} and \ref{thm:qc-englulf_2}
one obtains

\begin{cor} \label{cor:engulf+t-f-gferf} Let $G$ be a torsion-free
hyperbolic group where each proper quasiconvex subgroup is engulfed. Then $G$ is GFERF.
\end{cor}


N. Romanovski\v\i{}  \cite{Romanovskii} and, independently, R. Burns \cite{Burns} showed that a free
product of two LERF groups is again a LERF group. We suggest yet another way for constructing GFERF groups
by proving a similar result for them (see Section \ref{sec:free-prod-gferf}):

\begin{thm} \label{thm:free-prod-gferf} Suppose $G_1$ and $G_2$  are GFERF hyperbolic groups. Then the free product
$G=G_1*G_2$ is also a GFERF hyperbolic group.
\end{thm}

{\bf Acknowledgements}

I would like to thank Dr. Alexander Ol'shanskii for useful discussions and Dr. Graham Niblo
for referring me to his paper \cite{Niblo-Wil}.

\section{Preliminaries}
Let $G$ be a group with a finite symmetrized generating set $\cA$. Naturally, this generating set gives rise to a
{\it word length} function $|g|_G$ for every element $g \in G$. The (left-invariant)
{\it word metric} $d:G \times G \to \N \cup \{0\}$
is defined by the formula $d(x,y)=|x^{-1}y|_G$ for any $x,y \in G$. This metric can be canonically extended to
the Cayley graph {\ga} by making each edge isometric to the interval $[0,1]\subset \mathbb{R}$.

%

For any three points $x,y,w \in \Gamma(G,\mathcal{A})$, the {\it Gromov product} of $x$ and $y$ with respect to $w$
is defined as $$(x|y)_w \stackrel{def}{=} \frac12 \Bigl(d(x,w)+d(y,w)-d(x,y) \Bigr).$$

Since the metric is left-invariant, for arbitrary $x,y,w \in G$ we have
$$(x|y)_w=(w^{-1}x|w^{-1}y)_{1_G}.$$

The group $G$ is called ({\it word}) {\it hyperbolic} according to M. Gromov \cite{Gromov} if
there exists $\delta \ge 0$ such that for any $x,y,z,w \in \Gamma(G,\cA)$ their Gromov products satisfy
$$(x|y)_w \ge min\{(x|z)_w,(y|z)_w\} - \delta .$$
Equivalently, $G$ is {\it hyperbolic} if there exists $\delta \ge 0$ such that each geodesic triangle
$\Delta$ in {\ga} is $\delta$-slim, i.e., every side of $\Delta$ is contained an a $\delta$-neighborhood of
the two other sides (see \cite{Mihalik}).

From now on we will assume that $G$ is a hyperbolic group and the constant $\delta$ is large
enough so that it satisfies both of the definitions above.

For any two points $x,y \in \Gamma(G,{\mathcal A})$ we fix a geodesic path between them and denote it by $[x,y]$.
Let $p$ be a path in the Cayley graph of $G$. Then $p_-$, $p_{+}$  will denote the startpoint and
the endpoint of $p$, $\|p\|$ -- its length; $lab(p)$, as usual, will mean the word in the alphabet $\mathcal A$
written on $p$. $elem(p) \in G$ will denote the element of the group $G$ represented by the word $lab(p)$.
If $W$ is a word in the $\mathcal A$, $elem(W)$ will denote the corresponding element of the group $G$.
For a subset $A \subset \Gamma(G,\cA)$ its closed $\varepsilon$-neighborhood will be denoted by
$\mathcal{O}_\varepsilon(A)$.

The $\delta$-slimness of geodesic triangles  implies $2\delta$-slimness
of all geodesic quadrangles $abcd$ in {\ga}:
$$[a,b] \subset \mathcal{O}_{2\delta}\bigl([b,c]\cup [c,d] \cup [a,d]\bigr).$$

A path $q$ is called $(\lambda,c)$-{ \it quasigeodesic} if there exist $0<\lambda \le 1$, $c \ge 0$, such that
for any subpath $p$ of $q$ the inequality $\lambda \|p\| - c \le d(p_-,p_+)$ holds. A word $W$ is said to be
$(\lambda,c)$-quasigeodesic if some (equivalently, every) path $q$ in {\ga} labelled by $W$ is
$(\lambda,c)$-quasigeodesic.


\begin{lemma} \label{close} {\normalfont (\cite[5.6,5.11]{Ghys},\cite[3.3]{Mihalik})}
There is a constant $\nu=\nu(\delta,\lambda,c)$ such that for any
$(\lambda,c)$-quasigeodesic path $p$ in $\Gamma(G,{\mathcal A})$ and a geodesic
$q$ with $p_- = q_-$, $p_+ = q_+$, one has $p \subset {\mathcal O}_\nu(q)$ and $q \subset
{\mathcal O}_\nu(p)$.
\end{lemma}

\begin{lemma} \label{quadrangle}  {\normalfont (\cite[Lemma 4.1]{paper2})} Consider  a geodesic
quadrangle $x_1x_2x_3x_4$  in the Cayley graph {\ga} with
$d(x_2,x_3)>d(x_1,x_2)+d(x_3,x_4)$. Then there are points $u,v \in [x_2,x_3]$ such that
$d(x_2,u) \le d(x_1,x_2)$, $d(v,x_3) \le d(x_3,x_4)$ and  the geodesic subsegment
$[u,v]$ of $[x_2,x_3]$ lies $2\delta$-close to the side $[x_1,x_4]$.
\end{lemma}

If $x,g \in G$, we define $x^g =gxg^{-1}$. For a subset $A$ of the group $G$, $A^g = gAg^{-1}$, and  the notation
$A^G$ will be used to denote the subset $\{gag^{-1}~|~a \in A, g \in G\} \subset G$.

\begin{remark} \cite[Remark 2.2]{hyp} \label{rem:qc-shifts} Let $Q\subseteq G$ be quasiconvex, $g\in G$. Then
the subsets $gQ$, $Qg$ and $gQg^{-1}$ are quasiconvex as well.
\end{remark}

Thus, a conjugate of a quasiconvex subgroup in a hyperbolic group is quasiconvex itself. Another important property of
hyperbolic groups states that any cyclic subgroup is quasiconvex (see \cite{Mihalik}, for instance).

In this paper we will also use the concept of {\it Gromov boundary} of a hyperbolic group $G$ (for a detailed theory
the reader is referred to the corresponding chapters in \cite{Ghys} or \cite{Mihalik}).)
A sequence $(x_i)_{i\in \N}$ of elements of the group $G$ is said to be {\it converging to infinity} if
$$\lim_{i,j \to \infty} (x_i|x_j)_{1_G} = \infty.$$

Two sequences $(x_i)_{i\in \N},(y_j)_{i\in \N}$ converging to infinity are said to be equivalent
if $$\lim_{i \to \infty} (x_i|y_i)_{1_G} = \infty.$$
The points of the boundary $\partial G$ are identified with the equivalence classes of sequences
converging to infinity. It is easy to see that this definition does not depend on the choice
of a basepoint: instead of $1_G$ one could use any fixed point $p$ of {\ga} (see \cite{Mihalik}).
If $\alpha$ is the equivalence class of $(x_i)_{i\in \N}$, we will write
$\displaystyle \lim_{i \to \infty} x_i = \alpha$.

The space $\partial G$ can be topologized so that it becomes compact, Hausdorff
and metrizable (see \cite{Mihalik},\cite{Ghys}).

Left multiplication by an element of the group induces a homeomorphic action of $G$ on its boundary:
for any $g \in G$ and $\left[(x_i)_{i\in \N}\right] \in \partial G$ set
 $$g \circ \left[(x_i)_{i\in \N}\right]\stackrel{def}{=}
\left[(gx_i)_{i\in \N}\right] \in \partial G.$$

If $g \in G$ is an element of infinite order then the sequences $(g^i)_{i\in \N}$ and
$(g^{-i})_{i\in \N}$ converge to infinity and we will use the notation
$$\lim_{i \to \infty} g^i = g^\infty \in \partial G,~\lim_{i \to \infty} g^{-i} = g^{-\infty} \in \partial G.$$

For a subset $A \subseteq G$ the {\it limit set} $\Lambda (A)$ of $A$ is the collection of the points
$\alpha \in \partial G$ that are limits of sequences (converging to infinity) from $A$.

Let us recall an auxiliary binary relation defined between subsets of an arbitrary group $G$ in \cite{paper2}:
suppose $A,B \subseteq G$, we will write $A \preceq B$ if and only if there exist elements
$x_1,\dots,x_n \in G$ such that $$A \subseteq Bx_1 \cup Bx_2 \cup \dots \cup Bx_n.$$

It is not difficult to see that the relation "$\preceq$" is transitive and for any $g \in G$, $A \preceq B$ implies
$gA \preceq gB$.

\begin{lemma} \label{sgpcompar} {\rm (\cite[Lemma 2.1]{paper2})}  Let $A,B$ be subgroups of $G$. Then
$A \preceq B$ if and only if the index $|A:(A \cap B)|$ is finite.
\end{lemma}

The basic properties of limit sets are described in the following statement:

\begin{lemma} \label{limitsets} {\normalfont (\cite{K-S},\cite{Swenson},\cite[Lemma 6.2]{paper2})}
Suppose $A,B$ are arbitrary subsets of $G$, $g \in G$. Then

{\rm (a)} $\Lambda (A) = \emptyset$ if and only if $A$ is finite;

{\rm (b)} $\Lambda (A)$ is a closed subset of the boundary $\partial G$;

{\rm (c)} $\Lambda(A \cup B) = \Lambda (A) \cup \Lambda (B)$;

{\rm (d)} $\Lambda(Ag) = \Lambda(A)$, $g \circ \Lambda(A)=\Lambda(gA)$;

{\rm (e)} if $A \preceq B$ then $\Lambda(A) \subseteq \Lambda(B)$. 
\end{lemma}

%

The following property of limit sets of quasiconvex subgroups was first proved by E. Swenson:
\begin{lemma}{\rm (\cite[Thm. 8]{Swenson}, \cite[Lemma 9.1]{disser})}\label{lem:inter_limsets}
Let $A$, $B$ be quasiconvex subgroups of a hyperbolic group $G$. Then
$\Lambda(A) \cap \Lambda(B)=\Lambda(A \cap B)$ in $\partial G$.
\end{lemma}

\section{Auxiliary Facts}

\begin{lemma}\label{lem:main-conj} Assume $H$ is an $\eta$-quasiconvex subgroup of a $\delta$-hyperbolic group $G$, $X$ is a word
over $\mathcal A$ representing an element of infinite order in $G$, $0< \lambda \le 1$ and $c \ge 0$.
Let $\nu=\nu(\delta,\lambda,c)$ be the constant given by Lemma \ref{close}.
There exists a number $N=N(\delta,\eta,\nu,G) \in \N$ such that for any $m \in \N$ the following property holds.

If a word $W \equiv U_1X^nU_2$ is $(\lambda,c)$-quasigeodesic and satisfies
$\|U_1\|,\|U_2\| > (m+\nu+c)/\lambda$, $n \ge N$ and $elem(V_1WV_2) \in H$ for some words $V_1,V_2$ with
$\|V_1\|,\|V_2\|\le m$, then there exist $k \in \N$ and $a \in G$ such that $elem(X^k)\in H^a$ and
$|a|_G \le 2\delta+\nu+\eta$.
\end{lemma}

\begin{proof} Consider a path $q$ in {\ga} starting at $1_G$ and labelled by $V_1WV_2$. By our assumptions, $q_+ \in H$.
Let $p$ and $r$ be its $(\lambda,c)$-quasigeodesic subpaths with $lab(p) \equiv W$ and $lab(r)\equiv X^n$ respectively.
Choose an arbitrary {\it phase} vertex $u \in r$ such that the subpath of $r$ from $r_-$ to $u$ is labelled by some power of $X$.

According to Lemma \ref{close} there is $v \in [p_-,p_+]$ satisfying $d(u,v) \le \nu$. Using the
assumptions and the triangle inequality one can achieve
$$d(p_-,v) \ge d(p_-,u)-d(u,v) \ge \lambda \|U_1\| -c -\nu > m, \mbox{ and }$$
$$d(p_+,v) \ge d(p_+,u)-d(u,v) \ge \lambda \|U_2\| -c -\nu > m.$$

Hence, by Lemma \ref{quadrangle}, $v \in {\mathcal O}_{2\delta}([1_G,q_+])$.

Thus, $$u \in {\mathcal O}_{2\delta+\nu}([1_G,q_+]) \subset {\mathcal O}_{2\delta+\nu+\eta}(H),$$
i.e., there is an element $a=a(u) \in G$ such that $|a|_G \le 2\delta+\nu+\eta$
and  $u \in Ha$.

Now, since the alphabet $\mathcal A$ is finite, there are only finitely many elements in $G$ having length at most
($2\delta+\nu+\eta$). Hence, if $n$ is large enough, there should be two different phase vertices $u_1,u_2 \in r$
with $a(u_1)=a(u_2)=a$. By the construction, $$u_1^{-1}u_2=elem(X^k)\in a^{-1}Ha$$ for some
positive integer $k$ ($X^k$ is a label of the segment of $r$ from $u_1$ to $u_2$, provided these points are chosen
in a correct order). Q.e.d.
\end{proof}

\begin{lemma}\label{lem:gh_le_h} Assume $G$ is a hyperbolic group and $H \le G$ is a quasiconvex subgroup.
If $g \in G$ and $gH \preceq H$, then $H \preceq gH$.
\end{lemma}

\begin{proof} If $H$ is finite, the statement is trivial. Our assumptions immediately imply
$$g^{k-1} H \preceq g^{k-2}H \preceq \dots \preceq gH \preceq H$$ for all
$k \in \N$. Hence, \begin{equation} \label{eq:gk-1h} g^{k-1}H \preceq H.
\end{equation}

If $g$ has order $k \in \N$ in the group $G$, then to achieve the desired result
it is enough to multiply both sides of the above formula by $g$.

Thus, we can further assume that $H$ is infinite and the element $g$ has infinite order. Therefore, $H$ has at least one
limit point $\alpha \in \Lambda(H)$. Observe that \eqref{eq:gk-1h} implies
$$g^n \circ \Lambda(H)=\Lambda(g^n H) \subseteq \Lambda(H),$$ thus $g^n \circ \alpha \in \Lambda(H)$ for all $n \in \N$.

If $\alpha \neq g^{-\infty}$ in $\partial G$, it is well known (see, for instance, \cite[8.16]{Ghys}) that the sequence
$(g^n \circ \alpha)_{n\in \N}$ converges to $g^\infty$. Since $\Lambda(H)$ is a closed subset of the boundary
$\partial G$, in either case we achieve
$$\Lambda(\langle g \rangle_\infty) =\{g^\infty,g^{-\infty}\} \cap \Lambda(H) \neq \emptyset.$$
By Lemma \ref{lem:inter_limsets}, the latter implies $g^k \in H$ for some $k \in \N$. Combining this fact with
\eqref{eq:gk-1h} we get $H=g^kH \preceq gH$, which concludes the proof.
\end{proof}

The previous lemma has an interesting corollary:
\begin{lemma}\label{lem:comm-of-qc} Suppose $H,K$ are subgroups of a hyperbolic group $G$ such that  $H \le K$,
$|K:H|=\infty$ and $H$ is quasiconvex. Then $|K:(K \cap H^g)|=\infty$ for any $g \in G$.
\end{lemma}

\begin{proof} If $|K:(K \cap H^g)|<\infty$ for some $g \in G$ then
$H \preceq K \preceq gHg^{-1} \preceq gH$ according to Lemma \ref{sgpcompar}. Consequently $g^{-1}H \preceq H$. Hence
$H \preceq g^{-1}H$ by Lemma \ref{lem:gh_le_h}, implying $gHg^{-1} \preceq gH \preceq H$. But the latter leads to
$K \preceq H$ which contradicts the condition $|K:H|=\infty$ (see Lemma \ref{sgpcompar}).
\end{proof}

If $G$ is a hyperbolic group, each element $g \in G$ of infinite order
belongs to a unique {\it maximal elementary subgroup} $E(g)$.
By \cite[Lemmas 1.16, 1.17]{Olsh2} it has the following description:
\begin{multline} E(g) = \bigl\{x \in G~|~xg^kx^{-1}=g^{l}~\mbox{for some $k,l \in \Z\backslash \{0\}$}\bigr\}
= \\ \bigl\{x \in G~|~xg^nx^{-1}=g^{\pm n}~\mbox{for some } n \in \N\bigr\}.
\label{elemdef}
\end{multline}
Note that the subgroup $\displaystyle E^+(g)\stackrel{def}{=}\bigl\{~x \in G~|~xg^nx^{-1}=g^{n}~\mbox{for some } n \in \N~\bigr\}$
has index at most 2 in $E(g)$.

Let $W_1,W_2,\dots,W_l$ be words in $\mathcal A$ representing elements $g_1,g_2,\dots,g_l$
of infinite order, where $E(g_i) \neq E(g_j)$ for $i \neq j$. The following lemma
will be useful:

\begin{lemma} \label{quasigeodesic} {\normalfont (\cite[Lemma 2.3]{Olsh2})}
There exist constants
$\lambda = \lambda(W_1,W_2,\dots,W_l)>0$, $c=c(W_1,W_2,\dots,W_l) \ge 0$ and
$N = N(W_1,W_2,\dots,W_l)>0$ such that any path $p$ in the Cayley graph {\ga}
with label $W_{i_1}^{m_1}W_{i_2}^{m_2} \dots W_{i_s}^{m_s}$ is $(\lambda,c)$-quasigeodesic
if $i_k \neq i_{k+1}$ for $k=1,2,\dots,s-1$, and $|m_k|>N$ for $k=2,3,\dots,s-1$
(each $i_k$ belongs to $\{1,\dots,l\}$).
\end{lemma}

For a subgroup $H$ of $G$ denote by $H^0$ the set of elements of infinite order in $H$; if $A \subseteq G$,
the subgroup $C_H(A)=\{h \in H~|~g^h=g, ~\forall~ g \in A\}$ is the {\it centralizer} of $A$ in $H$.

Set $E(H)=\bigcap_{x \in H^0} E(x)$. If $H$ is a non-elementary subgroup of $G$, then $E(H)$ is the unique
maximal finite subgroup of $G$ normalized by $H$ (\cite[Prop. 1]{Olsh2}).

If $g \in G^0$, $T(g)$ will be used to denote the set of elements of finite order in the subgroup $E(g)$.

Let $G$ be a hyperbolic group and $H$ be its non-elementary subgroup.
Recalling the definition from \cite{Olsh2} (and using the terminology from \cite{paper3}), an element $g \in H^0$ will be called
$H$-{\it suitable} if $E(H)=T(g)$ and $$E(g)=E^+(g)=C_G(g)=T(g)\times\langle g \rangle.$$
In  particular, if the element $g$ is $H$-suitable then $g \in C_H\bigl(E(H)\bigr)$.

Two elements $g,h \in G$ of infinite order are called {\it commensurable} if  $g^k=ah^la^{-1}$ for some non-zero
integers $k,l$ and some $a \in G$. The following important statement was proved by A. Ol'shanskii in 1993:

\begin{lemma} {\rm (\cite[Lemma 3.8]{Olsh2})} \label{lemma3.8} Every non-elementary subgroup $H$ of a hyperbolic
group $G$ contains an infinite set of pairwise non-commensurable $H$-{\it suitable} elements.
\end{lemma}

Suitable elements can be modified in a natural way:

\begin{lemma} \label{y-mod} {\rm (\cite[Lemma 4.3]{paper3})} Let $H$ be a non-elementary subgroup of a hyperbolic
group $G$, and
$g$ be an $H$-suitable element. If $y \in C_H\bigl(E(H)\bigr) \backslash E(g)$ then there exists $N \in \N$
such that the element $yg^n$ has infinite order in $H$ and is $H$-suitable for every $n\ge N$.
\end{lemma}

In \cite{paper2} the author studied properties of quasiconvex subgroups of infinite index and showed

\begin{lemma} \label{lem:non-conj-elem} {\rm (\cite[Prop. 1]{paper2})}
Suppose $H$ is a quasiconvex subgroup of a hyperbolic group $G$ and $K$ is any subgroup of
$G$ that satisfies $|K:(K \cap H^g)| = \infty$ for all $g\in G$. Then there
exists an element $x \in K$ having infinite order, such that
$\langle x \rangle_\infty \cap H^G = \{1_G\}$.
\end{lemma}

Later we will utilize a stronger fact:

\begin{lemma}\label{lem:suit-mod-non-conj} Assume $H,K$ are subgroups of a $\delta$-hyperbolic group $G$,
$H$ is $\eta$-quasiconvex, $K$ is non-elementary and $|K:(K\cap H^g)|=\infty$ for every $g \in G$.
Then there exists a $K$-suitable element $y \in K$ such that $\langle y \rangle_\infty \cap H^G=\{1_G\}$.
\end{lemma}

\begin{proof} Set $K'=C_K \bigl( E(K) \bigr)$. 
Since $E(K)$ is a finite subgroup normalized by $K$, we have $|K:K'|<\infty$.
Therefore $|K':(K'\cap H^g)|=\infty$ for all $g \in G$. Applying
Lemma \ref{lem:non-conj-elem}, we can find an element of infinite order $x \in K'$ such that
$\langle x \rangle_\infty \cap H^G = \{1_G\}$. By Lemma \ref{lemma3.8} there is a $K$-suitable element $z\in K$
which is non-commensurable with $x$, hence $\langle x \rangle_\infty \cap E(z)=\{1_G\}$.

Choose some words $X$ and $Z$ in the alphabet $\mathcal A$ representing $x$ and $z$ respectively. Then one is able
to find the numbers $\lambda=\lambda(X,Z)$, $c=c(X,Z)$ and $N_1=N_1(X,Z)$ from the claim of Lemma \ref{quasigeodesic}.

Define $\nu=\nu(\delta,\lambda,c)$ and $N_2=N_2(\delta,\eta,\nu,G)$ as in Lemmas \ref{close} and \ref{lem:main-conj}.
Denote $N=\max \{N_1,N_2\}$ and apply Lemma \ref{y-mod} to obtain $n \ge N_1$ such that the element $y=x^Nz^n \in K$ is
$K$-suitable.

It remains to check that $\langle y \rangle_\infty \cap H^G=\{1_G\}$. Assume the contrary, i.e., there exist $t \in \N$
and $g\in G$ such that $y^t \in H^g$. Then for each $l \in \N$, the element $y^{tl} \in H^g$ will be represented by the
$(\lambda,c)$-quasigeodesic word $W \equiv (X^NZ^n)^{tl}$. And if the number $l$ is chosen sufficiently large
(compared to $m=|g^{-1}|_G=|g|_G$),
it should be possible to find a subword of the form $X^N$ in the "middle" of $W$ which
satisfies all the assumptions of Lemma \ref{lem:main-conj}. Hence, $x^k=elem(X^k)\in H^G$ for some $k \in \N$.
The latter contradicts to the construction of $x$. Thus the lemma is proved.
\end{proof}

As usual, let $G$ be a $\delta$-hyperbolic group and $H$-- its $\eta$-quasiconvex subgroup.

\begin{lemma} \label{lem:non-conj-free} Suppose the elements $x_1,x_2 \in G$ have infinite order, $E(x_1)\neq E(x_2)$,
and for each $i=1,2$, satisfy
$\langle x_i \rangle_\infty \cap H^G=\{1_G\}$. Then there exists a number $N \in \N$ such that for any $m,n \ge N$
the elements $x_1^m$, $x_2^n$ freely generate a free subgroup $F$ of rank $2$ in $G$. Moreover,
$F$ is quasiconvex and $F\cap H^G=\{1_G\}$.
\end{lemma}

\begin{proof} Choose some words $X_1$ and $X_2$ in the alphabet $\mathcal A$
representing the elements $x_1$ and $x_2$. Apply Lemma \ref{quasigeodesic} to find the corresponding
$\lambda=\lambda(X_1,X_2)$, $c=c(X_1,X_2)$ and $N_1=N_1(X_1,X_2)$. Then one can find the constant
$\nu=\nu(\delta,\lambda,c)$ from the claim of Lemma \ref{close} and define $N_2=N_2(\delta,\eta,\nu,G)$
in accordance with Lemma \ref{lem:main-conj}.

Set $N=\max\{N_1,N_2, \lfloor c/\lambda \rfloor+1\}$. Consider arbitrary integers $m,n \ge N$ and the subgroup
$F=\langle x_1^m,x_2^n\rangle \le G$.
By Lemma \ref{quasigeodesic} any non-empty (freely) reduced word $W$ in the generators $\{X_1^m, X_2^n\}$ is
$(\lambda,c)$-quasigeodesic.
Hence $$|elem(W)|_G \ge \lambda \|W\|-c \ge \lambda N -c>0.$$ Consequently, $elem(W) \neq 1_G$ in $G$, implying that
$F$ is free with the free generating set $\{x_1^m,x_2^n\}$. By the construction of $\nu$, $F$ will be
$\varepsilon$-quasiconvex where $$\varepsilon=\nu+ \frac 12 \max\{m \|X_1\|,n \|X_2\|\}.$$

Consider a non-empty cyclically reduced word $W$ in the generators $\{X_1^m, X_2^n\}$.
For establishing the last claim, it is
sufficient to demonstrate that $elem(W) \notin H^G$. Arguing as in Lemma \ref{lem:suit-mod-non-conj}, suppose
$elem(W) \in H^g$ for some $g\in G$. Then $\left(elem(W)\right)^l \in H^g$ for every $l \in \N$. By choosing $l$
sufficiently large and applying Lemma \ref{lem:main-conj} one will obtain a contradiction with the assumption
$\langle x_i \rangle_\infty \cap H^G=\{1_G\}$, $i=1,2$, as in Lemma \ref{lem:suit-mod-non-conj}.
Therefore $F \cap H^G=\{1_G\}$. Q.e.d.
\end{proof}

\begin{cor} \label{cor:free-non-conj-good}
With the assumptions of Lemma \ref{lem:suit-mod-non-conj}, $K$ has a free subgroup $F$ of rank $2$
which is quasiconvex in $G$, $E(F)=E(K)$ and $F \cap H^G=\{1_G\}$.
\end{cor}

\begin{proof} Choose a $K$-suitable element $x_1 \in K$ according to Lemma \ref{lem:suit-mod-non-conj}. Since $K$ is
non-elementary, there exists $y \in K \backslash E(x_1)$. Therefore, $x_2 \stackrel{def}{=}yx_1y^{-1} \in K^0$ and
$E(x_2)\neq E(x_1)$ (as it can be seen from \eqref{elemdef}). By the construction,
$\langle x_i \rangle_\infty \cap H^G=\{1_G\}$, $i=1,2$. Hence the subgroup $F$ can be found by applying
Lemma \ref{lem:non-conj-free}. Evidently $E(K) \le E(F)$, and $E(F)\subseteq T(x_1)=E(K)$. Thus $E(F)=E(K)$.
\end{proof}

\section{Free products of quasiconvex subgroups}
Below we will assume that $G$ is a $\delta$-hyperbolic group generated by
a finite set $\mathcal A$. 

First let us recall some properties of the hyperbolic boundary.

\begin{lemma} \label{boundedproduct} {\rm (\cite[Lemma 2.14]{paper3})}
Suppose $A$ and $B$ are arbitrary subsets of $G$ and
$\Lambda(A)\cap \Lambda(B)=\emptyset$. Then $\displaystyle \sup_{a \in A,~b\in B}\{(a|b)_{1_G}\}< \infty$.
\end{lemma}

%

\begin{remark} \label{rem:inf-ord-elem-2} Suppose $g,x \in G$ and $g$ has infinite order.
If $\left( x \circ \{g^{\pm \infty}\} \right) \cap  \{g^{\pm \infty}\}  \neq \emptyset $ in the boundary
 $\partial G$, then $x \in E(g)$. If $E(g)=E^+(g)$ then $g^\infty \notin G \circ \{g^{-\infty}\}$.
\end{remark}

Note that  $x \circ \{g^{\pm \infty}\}=\{(xgx^{-1})^{\pm \infty}\}=\Lambda(\langle xgx^{-1} \rangle)
\subset \partial G$. Since any cyclic subgroup in a hyperbolic group is quasiconvex, we can apply
Lemma \ref{lem:inter_limsets} to show that $\langle g \rangle \cap \langle xgx^{-1} \rangle \neq \{1_G\}$.
Hence $x \in E(g)$. \\
Let $E(g)=E^+(g)$ and suppose that $g^\infty = x \circ g^{-\infty}$ for some $x \in G$.
Then, as we showed above, $x \in E(g)=E^+(g)$. Hence
$$x \circ g^{-\infty}=\lim_{n \to -\infty} (xgx^{-1})^n=\lim_{n \to -\infty} g^n=g^{-\infty}.$$
Thus we achieve a contradiction with the inequality $g^\infty \neq g^{-\infty}$.

\begin{lemma} \label{lem:mod-inf-ord} Let $g,x \in G$ where $g$ has
infinite order and $E(g)=E^+(g)$. Then there is $N_1 \in \N$ such that for every $n \ge N_1$ the element $xg^n \in G$ has infinite order.
\end{lemma}

\begin{proof} If $x \notin E(g)$ then the claim follows by
\cite[Lemma 9.14]{disser}.

So, assume $x \in E(g)$. Since $E(g)=E^+(g)$,
the center of $E(g)$ has finite index in it, thus $E(g)$ is an FC-group.
By B.H. Neumann's Theorem \cite{Neumann-FC} the elements of finite order form a subgroup $T(g) \le E(g)$. Therefore
the cardinality of the intersection $\{xg^k~|~k\in \Z\} \cap T(g)$ can be at most $1$. Thus $xg^n \notin T(g)$ for each sufficiently large $n$.
\end{proof}

The main result of this paper is based on the following statement concerning broken lines in a
$\delta$-hyperbolic metric space:

\begin{lemma} \label{brokenlines} {\normalfont (\cite[Lemma 21]{Olsh1},
\cite[Lemma 3.5]{paper2}) } Let $p = [y_0,y_1,\dots,y_n]$ be a broken
line in $\Gamma(G,{\mathcal A})$ such that
$||[y_{i-1},y_i]|| > C_1$ $\forall~i=1,\dots,n$, and $(y_{i-1}|y_{i+1})_{y_i} < C_0$ $\forall~i=1,\dots,n-1$,
where $C_0 \ge 14 \delta$, $C_1 > 12(C_0 + \delta)$. Then the geodesic segment $[y_0,y_n]$ is contained in the closed $14\delta$-neighborhood
of $p$ and $\|[y_0,y_n]\| \ge \|p\|/2$.
\end{lemma}

Suppose $a,b,c$ and $d$ are arbitrary points in {\ga}. Considering the geodesic triangle with the vertices
$1_G, a$ and $ab$, one can observe that
$$(a|ab)_{1_G}=|a|_G-(1_G|ab)_a=|a|_G-(a^{-1}|b)_{1_G}.$$
Now, since {\ga} is $\delta$-hyperbolic,
\begin{equation}\label{eq:ineq1}
(a|c)_{1_G} \ge \min\bigl\{(a|ab)_{1_G},(ab|c)_{1_G}\bigr\}-\delta =
\min\bigl\{|a|_G-(a^{-1}|b)_{1_G},(ab|c)_{1_G}\bigr\}-\delta.
\end{equation}

Replacing $c$ with $cd$ in the above formula, we get
\begin{equation}\label{eq:ineq2}
(a|cd)_{1_G} \ge \min\bigl\{|a|_G-(a^{-1}|b)_{1_G},(ab|cd)_{1_G}\bigr\}-\delta.
\end{equation}

The Gromov product is symmetric, therefore we are able to combine formulas
\eqref{eq:ineq1} and \eqref{eq:ineq2} to achieve
\begin{multline}\label{eq:ineq3}
(a|c)_{1_G} \ge \min\bigl\{|c|_G-(c^{-1}|d)_{1_G},(a|cd)_{1_G}\bigr\}-\delta
\ge \\ \min\bigl\{|a|_G-(a^{-1}|b)_{1_G},|c|_G-(c^{-1}|d)_{1_G}, (ab|cd)_{1_G}\bigr\}-2\delta.
\end{multline}

\begin{thm} \label{thm:main-engulf} Consider some elements $g_1,x_1,g_2,x_2,\dots,g_s,x_s \in G$ and
an $\eta$-quasi\-convex subgroup $H \le G$.
Suppose the following three conditions are satisfied:
\begin{itemize}
\item $g_1,\dots, g_s$ have infinite order and are pairwise
non-commensurable; \item $E(g_i)=E^+(g_i)$ for each
$i=1,2,\dots,s$; \item $E(g_i) \cap H = E(g_i) \cap x_i^{-1} H x_i =\{1_G\}$ for each
$i=1,2,\dots,s$.
\end{itemize}
Then there exists a number $N \in \N$ such that for every $n \ge N$ the
elements $x_ig_i^n \in G$, $i=1,2\dots,s$, have infinite order, and the
subgroup $$M \stackrel{def}{=}\langle H, x_1g_1^n,\dots,x_sg_s^n \rangle \le G$$ is quasiconvex in $G$ and
isomorphic (in the canonical way) to the free product
$H * \langle x_1g_1^n \rangle * \dots * \langle x_sg_s^n \rangle$.
\end{thm}

\begin{proof} Choose arbitrary elements $w_1,w_2 \in M$ and define
$w=w_1^{-1}w_2 \in M$. Then this element has a presentation
\begin{equation} \label{eq:w}
w=h_0 (x_{i_1}g_{i_1}^n)^{\epsilon_1} h_1 (x_{i_2}g_{i_2}^n)^{\epsilon_2} \cdots h_{l-1}
(x_{i_l}g_{i_l}^n)^{\epsilon_l}h_l,
\end{equation}
where $h_j \in H$, $i_j \in \{1,\dots,s\}$, $\epsilon_j \in \{1,-1\}$,
$j=1,2,\dots,l$, $l \in \N \cup \{0\}$.

Moreover, we can assume that the presentation \eqref{eq:w} is reduced in the
following sense: if $1 \le j \le l-1$, $i_j=i_{j+1}$ and
$\epsilon_{j+1}=-\epsilon_j$ then $h_j\neq 1_G$.

Consider a geodesic  broken line $[y_0,y_1,\dots, y_{l+1}]$ in {\ga} with $y_0=w_1$
 and $elem([y_k,y_{k+1}])=h_k (x_{i_{k+1}}g_{i_{k+1}}^n)^{\epsilon_{k+1}}$,
$k=0,1,\dots,l-1$, $elem([y_l,y_{l+1}])=h_l$. Therefore
$elem([y_0,y_{l+1}])=w$ and $y_{l+1}=w_2$.

Now we are going to find upper bounds for the Gromov products
$$(y_{k-1}|y_{k+1})_{y_k}=(y_k^{-1}y_{k-1}|y_k^{-1}y_{k+1})_{1_G}, ~k=1,\dots,l.$$

By the assumptions of the theorem, Lemma \ref{lem:inter_limsets} implies that
 $$x_i \circ g_i^\infty=(x_i g_i x_i^{-1})^\infty \notin \Lambda(H)~\mbox{ and }~g_i^{-\infty} \notin \Lambda(H),~i=1,\dots,s.$$

Since $$\Lambda \left( \{x_ig_i^m,g_i^{-m}x_i^{-1}~|~m \in \N, 1\le i \le s\} \right)=
\{x_i \circ g_i^\infty,g_i^{-\infty}~|~1\le i \le s\} \subset \partial G,$$
we are able to apply 
Lemma \ref{boundedproduct} to define
$$\alpha\stackrel{def}{=}\max\left\{(h|x_ig_i^m)_{1_G},(h|g_i^{-m}x_i^{-1})_{1_G}
~\bigl|~ h \in H, 1\le i \le s, m \in \N \right\}< \infty.$$

Similarly, since $g_i$ and $g_j$ are non-commensurable if $i \neq j$ and $E(g_i)=E^+(g_i)$, we have
(according to Lemma \ref{lem:inter_limsets} and Remark \ref{rem:inf-ord-elem-2})
$$G \circ\{g_i^{\pm \infty} \}\cap G\circ \{g_j^{\pm \infty}\}=\emptyset,
~G \circ\{g_i^{\infty} \}\cap G\circ \{g_i^{- \infty}\}=\emptyset,
~i,j \in \{1,\dots,s\}, i \neq j.$$
Hence, the following numbers are also finite: \begin{multline*} \beta_1\stackrel{def}{=}
\max\Bigl\{\bigl((x_ig_i^m)^{-1}|hx_jg_j^t \bigr)_{1_G}, ~\bigl(x_ig_i^m|h(x_jg_j^t)^{-1}
\bigr)_{1_G} \Bigl| \\ h\in H, ~|h|_G \le 2\alpha+2\delta,~1\le i,j \le s,
~m,t \in \N \Bigr\},\end{multline*}
\begin{multline*} \beta_2\stackrel{def}{=}
\max\Bigl\{\bigl(x_ig_i^m|hx_jg_j^t \bigr)_{1_G}, ~\bigl((x_ig_i^m)^{-1}|h(x_jg_j^t)^{-1}
\bigr)_{1_G} \Bigl| \\ h\in H, ~|h|_G \le 2\alpha+2\delta, ~1\le i,j \le s,
~i \neq j, ~m,t \in \N \Bigr\}.\end{multline*}

Note that if $h \in H \backslash \{1_G\}$ then, according to our assumptions,
$x_i^{-1}hx_i \notin E(g_i)$. Therefore,
$\{g_i^{\pm \infty}\} \cap (x_i^{-1}hx_i) \circ \{g_i^{\pm \infty}\}=\emptyset$ (Remark \ref{rem:inf-ord-elem-2}), implying
$x_i \circ g_i^\infty \neq (hx_i) \circ g_i^\infty$, $i=1,\dots,s$.
We can also use Remark \ref{rem:inf-ord-elem-2} to show that $g_i^{-\infty} \neq h \circ g_i^{-\infty}$ for each $i$.
Consequently, by Lemma \ref{boundedproduct},
\begin{multline*} \beta_3\stackrel{def}{=}
\max\Bigl\{\bigl(x_ig_i^m|hx_ig_i^t \bigr)_{1_G}, ~\bigl((x_ig_i^m)^{-1}|h(x_ig_i^t)^{-1}
\bigr)_{1_G} \Bigl| \\ h\in H,~ |h|_G \le 2\alpha+2\delta,~ h \neq 1_G,
~1\le i \le s,~m,t \in \N \Bigr\} < \infty.\end{multline*}

Finally, define $\beta=\max\{\beta_1,\beta_2,\beta_3\}<\infty,$
\begin{equation} \label{eq:C_0,C_1}
C_0=\max\{\alpha+2\delta,\beta+\delta,14\delta\}+1 ~\mbox{ and }~ C_1=12(C_0+\delta)+1.
\end{equation}

Since each $g_i$, $i=1,\dots,s$, has infinite order in $G$ there exists
$N \in \N$ such that for any $i \in \{1,\dots,s\}$ one has
\begin{equation} \label{eq:N}
|g_i^n|_G >  \max\{\alpha,\beta,2C_1\}+\alpha+|x_i|_G+2\delta,~\forall~n \ge N.
\end{equation}

By Lemma \ref{lem:mod-inf-ord}, without loss of generality, we can assume that the order of $x_ig_i^n$, $i=1,2,\dots,s$,
is infinite for every $n \ge N$.
Fix an integer $n \ge N$ and choose any $k \in \{1,\dots, l-1\}$. Then
 $$(y_{k-1}|y_{k+1})_{y_k}=
\bigl( (x_{i_k}g_{i_k}^n)^{-\epsilon_k} h_{k-1}^{-1} \bigr|
h_k (x_{i_{k+1}}g_{i_{k+1}}^n)^{\epsilon_{k+1}}\bigr)_{1_G}.$$

To simplify the notation, let us denote $a=(x_{i_k}g_{i_k}^n)^{-\epsilon_k}$,
$b=h_{k-1}^{-1}$, $c=h_k$ and $d=(x_{i_{k+1}}g_{i_{k+1}}^n)^{\epsilon_{k+1}}$.
By construction
\begin{equation} \label{eq:(a|c)etal} (a|c)_{1_G},(a^{-1}|b)_{1_G},(c^{-1}|d)_{1_G} \le \alpha. \end{equation}

Now we need to consider two separate cases.

{\bf Case 1:} $|h_k|_G=|c|_G \le 2\alpha+2\delta$. Then, due to the definition of the number  $\beta$,
the inequality $(a|cd)_{1_G} \le \beta$ holds.
Therefore, applying formulas \eqref{eq:ineq2} and \eqref{eq:(a|c)etal} one obtains
\begin{multline*} \beta \ge (a|cd)_{1_G} \ge \min \left\{ |a|_G-(a^{-1}|b)_{1_G},(ab|cd)_{1_G}\right\} -\delta \ge \\
\min \left\{ |g_{i_k}^n|_G-|x_{i_k}|_G-\alpha,(ab|cd)_{1_G}\right\}-\delta .\end{multline*}

By \eqref{eq:N}, $|g_{i_k}^n|_G-|x_{i_k}|_G-\alpha> \beta+\delta$, hence the above inequality implies
$$(y_{k-1}|y_{k+1})_{y_k}=(ab|cd)_{1_G} \le \beta + \delta < C_0.$$

{\bf Case 2:} $|h_k|_G=|c|_G > 2\alpha+2\delta$. Apply formulas \eqref{eq:ineq3} and \eqref{eq:(a|c)etal} to achieve
$$\alpha \ge (a|c)_{1_G} \ge \min \left\{ |a|_G-\alpha,|c|_G-\alpha,(ab|cd)_{1_G}\right\} -2\delta.$$

Observing  $|a|_G-\alpha\ge |g_{i_k}^n|_G-|x_{i_k}|_G-\alpha > \alpha+2\delta$ and
$|c|_G-\alpha > \alpha+2\delta$, we can conclude that
$$(y_{k-1}|y_{k+1})_{y_k}=(ab|cd)_{1_G} \le \alpha + 2\delta < C_0.$$

At last, let us estimate the product $(y_{l-1}|y_{l+1})_{y_l}=\bigl( (x_{i_l}g_{i_l}^n)^{-\epsilon_k}
h_{l-1}^{-1} \bigr| h_l \bigr)_{1_G}.$ Denote $a=(x_{i_l}g_{i_l}^n)^{-\epsilon_k}$, $b=h_{l-1}^{-1}$ and
$c=h_l$.

Using formula \eqref{eq:ineq1} and the definition of $\alpha$ one can obtain
$$\alpha \ge (a|c)_{1_G} \ge \min \left\{ |x_{i_l}g_{i_l}^n|_G-\alpha,(ab|c)_{1_G}\right\} -\delta.$$
As before, the latter implies that $$(y_{l-1}|y_{l+1})_{y_l} = (ab|c)_{1_G} \le \alpha+\delta< C_0.$$

Thus, we have shown that \begin{equation}\label{eq:angle-ineq}
(y_{k-1}|y_{k+1})_{y_k} < C_0 ~\mbox{for each } k=1,2,\dots,l.
\end{equation}

Now we need to find a lower bound for the lengths of the sides in the broken line
$[y_0,y_1,\dots,y_{l+1}]$.

Let $0 \le k \le l-1$. Note that $|ab|_G \ge |b|_G-(a^{-1}|b)_{1_G}$ for any $a,b \in G$. Hence
\begin{multline*} \|[y_k,y_{k+1}]\|=|h_k (x_{i_{k+1}}g_{i_{k+1}}^n)^{\epsilon_{k+1}}|_G \ge \\
|(x_{i_{k+1}}g_{i_{k+1}}^n)^{\epsilon_{k+1}}|_G-\bigl(h_k^{-1} \bigl|
(x_{i_{k+1}}g_{i_{k+1}}^n)^{\epsilon_{k+1}} \bigr)_{1_G} \ge |g_{i_{k+1}}^n|_G-|x_{i_{k+1}}|_G-\alpha . \end{multline*}

Applying inequality \eqref{eq:N} we obtain
\begin{equation} \label{eq:length-bound} \|[y_k,y_{k+1}]\| > 2C_1 ~\mbox{ if } 0 \le k \le l-1.
\end{equation}

Depending on the length of the last side, $\|[y_l,y_{l+1}]\|=|h_l|_G$, there can occur two
different situations.

{\bf Case 1:} $\|[y_l,y_{l+1}]\|=|h_l|_G \le C_1$. Then we can use inequalities \eqref{eq:angle-ineq} and
\eqref{eq:length-bound} to apply Lemma \ref{brokenlines} to the geodesic broken line
$p'=[y_0,\dots,y_l]$. Hence $[y_0,y_l]\subset {\mathcal O}_{14\delta}(p')$ and
$d(y_0,y_l) \ge \|p'\|/2$.

Since geodesic triangles in {\ga} are $\delta$-slim, we have
$$[y_0,y_{l+1}]\subset {\mathcal O}_{\delta}\bigl([y_0,y_l]\cup [y_l,y_{l+1}]\bigr) \subset
{\mathcal O}_{\delta+C_1}\bigl([y_0,y_l]\bigr) \subset {\mathcal O}_{15\delta+C_1}(p').$$

Now, if $l \ge 1$ in the presentation \eqref{eq:w}, one can use \eqref{eq:length-bound} to obtain
$$d(y_0,y_{l+1}) \ge d(y_0,y_l)-d(y_l,y_{l+1})\ge \|p'\|/2-C_1 \ge \|[y_0,y_1]\|/2-C_1>0.$$

{\bf Case 2:} $\|[y_l,y_{l+1}]\|=|h_l|_G > C_1$. Then we can apply Lemma \ref{brokenlines} to the broken
line $p=[y_0,\dots,y_l,y_{l+1}]$, thus achieving
$$ [y_0,y_{l+1}] \subset {\mathcal O}_{14\delta}(p).$$ As before, if $l \ge 1$ one has
$$d(y_0,y_{l+1})\ge \|p\|/2 >0.$$

Thus, in either case we have established the following properties:
\begin{equation} \label{eq:close} [y_0,y_{l+1}] \subset {\mathcal O}_{15\delta+C_1}(p) ~\mbox{ and }
\end{equation}
\begin{equation} \label{eq:long} d(y_0,y_{l+1}) >0.
\end{equation}

The inequality \eqref{eq:long} implies that $w \neq 1_G$ in the group $G$ for any element
$w \in M$ having a "reduced" presentation \eqref{eq:w} with $l \ge 1$. Therefore
$$M \cong H * \langle x_1g_1^n \rangle * \dots * \langle x_sg_s^n \rangle .$$

As $n \ge N$ is fixed, one is able to define the constants
$$\zeta=\max_{1\le i \le s} \{|x_ig_i^n|_G\}<\infty ~\mbox{ and }~
\varepsilon=16\delta+C_1+\eta+\zeta.$$
We will finish the proof by showing that $[y_0,y_{l+1}] \subset {\mathcal O}_{\varepsilon}(M)$
which implies that $M$ is $\varepsilon$-quasiconvex.

Consider an arbitrary $k \in \{0,1,\dots,l\}$.
Using the construction of $y_k \in M$ and $\zeta$, and $\delta$-hyperbolicity of the Cayley graph
we achieve
\begin{equation} \label{eq:y_k-h_k}
[y_k,y_{k+1}] \subset {\mathcal O}_{\delta+\zeta} ([y_k,y_kh_k]).
\end{equation}
$H$ is $\eta$-quasiconvex, therefore $[1_G,h_k] \subset  {\mathcal O}_{\eta}(H)$.
The metric in {\ga} is invariant under the action of $G$ by left translations, consequently
$$[y_k,y_kh_k] \subset  {\mathcal O}_{\eta}(y_kH)\subset  {\mathcal O}_{\eta}(M).$$

Combining the latter formula with \eqref{eq:y_k-h_k} leads us to
$$[y_k,y_{k+1}]\subset {\mathcal O}_{\delta+\eta+\zeta}(M)~ \mbox{ for each } k.$$

Finally, an application of \eqref{eq:close} yields
$$[y_0,y_{l+1}] \subset {\mathcal O}_{15\delta+C_1+\delta+\eta+\zeta}(M)=
{\mathcal O}_{\varepsilon}(M),$$
as desired.
\end{proof}

\section{Hyperbolic groups with engulfing}\label{sec:engulf}

\begin{proof}[Proof of Theorem \ref{thm:free-engulf}.] Since every elementary group is residually finite, it is
sufficient to consider the
case when $G$  is non-elementary. Let $x_1,\dots,x_s$ be the generators of $G$.

Define $\displaystyle K=\bigcap_{L \le G, |G:L|<\infty} L$; then $K$ is normal in $G$.
Suppose $K$ is infinite. If the subgroup $K$ were elementary, then it would be quasiconvex
(this follows directly from Lemmas \ref{quasigeodesic} and \ref{close}). Hence, according to a
result proved by Mihalik and Towle \cite{Towle} (see also \cite[Cor. 2]{paper2}), it would have a finite index
in $G$, thus the group $G$ would also have to be elementary. Therefore $K$ can not be elementary.

Now we can apply Lemma \ref{lemma3.8} to find pairwise non-commensurable  $K$-suitable elements
$g_1,\dots,g_{s+1} \in K$. Since the trivial subgroup $H=\{1_G\} \le G$ is quasiconvex, we can use
Theorem \ref{thm:main-engulf} to show that the subgroups
$M=\langle x_1g_1^n,\dots, x_sg_s^n \rangle$ and $M'=\langle x_1g_1^n,\dots, x_sg_s^n,g_{s+1}^n \rangle$
are free (of ranks $s$ and $(s+1)$ respectively) and quasiconvex in $G$ for some sufficiently large $n \in \N$.

Note that $M$ is a proper (infinite index) subgroup of $G$ because $|M':M|=\infty$. According to our
assumptions,  there exists a proper finite index subgroup $L$ in $G$ with $M \le L$. By the construction,
$K \le L$, thus $x_ig_i^n, g_i \in L$ for each $i=1,\dots,s$. Consequently $x_i \in L$ for each
$i=1,\dots,s$, contradicting with properness of $L$.

Therefore $K$ is finite. If $G$ is torsion-free then $K$ is trivial, and, thus, $G$ is residually finite.
\end{proof}

As $E(G)$ is the maximal finite normal subgroup in $G$, we obtain the following statement
right away:

\begin{cor} \label{cor:engulf-rf} With the assumptions of Theorem \ref{thm:free-engulf}, suppose, in addition,
that $E(G)=\{1_G\}$. Then $G$ is residually finite.
\end{cor}

Below it will be convenient to use the following equivalence relation between subsets of a group $G$
defined in \cite{paper2}: for any $A,B \subseteq G$ such that $A \preceq B$ and
$B \preceq A$ we will write $A \approx B$.

\begin{remark} \label{rem:equiv-qc} {\rm (\cite[Remark 3]{paper2})} If $A,B \subseteq G$,
$A \approx B$ and $A$ is quasiconvex, then $B$ is also quasiconvex.
\end{remark}

In particular, if $A \le B$ are subgroups of $G$ and $A$ has finite index in $B$ then $A \approx B$. Hence $A$
is quasiconvex if and only if $B$ is quasiconvex.

\begin{lemma} \label{lem:qc-englulf_1} Let $G$ be a residually finite hyperbolic group
and let $H \le G$ be a quasiconvex subgroup. Suppose that every proper quasiconvex subgroup of $G$ is engulfed.
Then $H$ has a finite index in its profinite closure $K$ in $G$.
\end{lemma}

\begin{proof} Arguing by the contrary, assume that $|K:H|=\infty$. Since
the group $G$ is residually finite, any finite subset is closed in the profinite topology.
Thus $H$ is infinite; hence $K$ is non-elementary.

Choose some generating set $x_1,\dots,x_s$ of the group $G$. Since $G$ is residually finite, it has a
finite index subgroup $G_1$ satisfying
\begin{equation} \label{eq:G_1} G_1 \cap \left( E(K) \cup \bigcup_{i=1}^s x_iE(K)x_i^{-1} \right)=\{1_G\}.
\end{equation}
Define $H_1 =H \cap G_1$; then $|H:H_1|<\infty$ and, according to Remark \ref{rem:equiv-qc},
$H_1$ is quasiconvex. The profinite closure $K_1$ of $H_1$ in $G$ has a finite
index in $K$, therefore $K_1$ is non-elementary and $|K_1:H_1|=\infty$. The definition of a finite index subgroup
implies that there is $l \in \N$ such that $g^l \in K_1$ for each $g \in K$. Since $E(g)=E(g^l)$ for any $g \in K^0$
we have
$$ E(K) \subseteq E(K_1) =
\bigcap_{g \in (K_1)^0} E(g) \subseteq \bigcap_{g \in K^0} E(g^l)= \bigcap_{g \in K^0} E(g)=E(K).$$
Thus $E(K_1)=E(K)$.

By Lemma \ref{lem:comm-of-qc}, $|K_1:(K_1 \cap H_1^g)|=\infty$ for every $g \in G$ and we can use
Corollary~\ref{cor:free-non-conj-good} to obtain a free subgroup $F \le K_1$ of rank $2$ satisfying
$F \cap H_1^G=\{1_G\}$ and $E(F)=E(K_1)=E(K)$.

 According to Lemma~\ref{lemma3.8} there
exist elements $g_1,\dots,g_{s+1} \in F$ which are pairwise non-commensurable and $F$-suitable.
Consequently $E(g_i)= \langle g_i \rangle_\infty \times E(K)$ and, since $E(K)$ is finite,
 formula \eqref{eq:G_1} implies
that $$E(g_i) \cap H_1= E(K) \cap H_1=\{1_G\},~i=1,2,\dots,s+1, ~\mbox{ and}$$
$$E(g_i) \cap x_i^{-1}H_1x_i=E(K) \cap x_i^{-1}H_1x_i=\{1_G\},~i=1,2,\dots,s.$$

Now we apply Theorem \ref{thm:main-engulf} to find $n \in \N$ such that the subgroups
$$M=\langle H_1, x_1g_1^n,\dots, x_sg_s^n \rangle ~\mbox{ and }~
M'=\langle H_1, x_1g_1^n,\dots, x_sg_s^n,g_{s+1}^n \rangle$$ are quasiconvex in $G$ and
$M'= M*\langle g_{s+1}^n\rangle_\infty \le G$. Thus, $|M':M|=\infty$ and $M$ is a proper
subgroup of $G$.

The subgroup $M$ is engulfed by our assumptions, therefore $G$ has a proper finite index subgroup
$L$ containing $M$. Observe that $K_1 \le L$ because $H_1 \le M \le L$, and, since
$x_ig_i^n,g_i \in M \cup K_1 \subset L$,
we get $x_i \in L$ for each $i=1,\dots,s$. The latter implies $G=L$ -- a contradiction.
\end{proof}

We will now prove Theorem \ref{thm:qc-englulf_2} which strengthens the the statement of the previous
lemma.

\begin{proof}[Proof of Theorem \ref{thm:qc-englulf_2}.] We can assume that $G$ is non-elementary because any
elementary group is LERF.
Since $G$ is residually finite, there is a finite index subgroup $G_1 \le G$ with
$G_1 \cap E(G) = \{1_G\}$.

Take an arbitrary quasiconvex subgroup $H \le G$ and set $H_1 = H \cap G_1$. As it follows from
Remark \ref{rem:equiv-qc}, $H_1$ is quasiconvex in $G$. Therefore, according to
Lemma \ref{lem:qc-englulf_1}, $H_1$ has a finite index in its profinite closure $K_1$
in $G$.

If $H_1=K_1$, i.e., $H_1$ is closed in the profinite topology on $G$, then so is $H$. Thus, we can suppose
that $H_1 \neq K_1$. Consequently $|G:H_1|=\infty$, and thus $|G:K_1|=\infty$.

The subgroup $K_1$ is quasiconvex according to Remark \ref{rem:equiv-qc}, hence we can apply
Lemma \ref{lem:suit-mod-non-conj} to find a $G$-suitable element $g \in G$ such that
$\langle g \rangle_\infty \cap K_1=\{1_G\}$. Since $E(g)=\langle g \rangle \times E(G)$, $E(G)$ is finite
and $K_1 \le G_1$, we have $$ E(g) \cap K_1=E(G) \cap K_1=\{1_G\}.$$

Now we can apply Theorem \ref{thm:main-engulf} to find a number $n\in \N$ such that the subgroups
$M=\langle H_1,g^n \rangle$ and $M'=\langle K_1,g^n \rangle$ are quasiconvex in $G$ and
$M' \cong K_1 * \langle g^n \rangle_\infty$.

Using properties of free products, we observe that  $M \le M'$ and $|M':M|=\infty$ because $H_1 \lneqq K_1$.
On the other hand, $M'$ is contained inside of the profinite closure of $M$ in $G$.
Thus we achieve a contradiction with the claim of Lemma \ref{lem:qc-englulf_1}.
\end{proof}

Before proceeding with the next statement, we need to recall some facts concerning quasiisometries
of metric spaces. Let $\mathcal X$ and $\mathcal Y$ be geodesic metric spaces with metrics $d(\cdot,\cdot)$ and
$e(\cdot,\cdot)$ respectively. A map $f: \mathcal{X} \to \mathcal{Y}$
is called a {\it quasiisometry} if there are constants $D_1 > 0$ and $D_2 \ge 0$ such that
$$ D_1^{-1}d(a,b)-D_2 \le e(f(a),f(b)) \le D_1 d(a,b)+D_2 \quad \forall~a,b \in \mathcal{X}.$$
The spaces $\mathcal{X}$ and $\mathcal Y$ are said to be {\it quasiisometric} if there exists
a quasiisometry $f:\mathcal{X} \to \mathcal{Y}$ whose image is {\it quasidense} in $\mathcal Y$, i.e.,
there exists $\varepsilon \ge 0$ such that for each $y \in \mathcal{Y}$ there is $x \in \mathcal{X}$ with
$e(y,f(x)) \le \varepsilon$.

M. Gromov \cite{Gromov} showed that if $\mathcal X$ is hyperbolic and quasiisometric to $\mathcal Y$
(through some map $f:\mathcal{X} \to \mathcal{Y}$) then the space $\mathcal Y$ is hyperbolic too.
He also noted that in this case the image $f(Q)$ of any quasiconvex subset $Q \subseteq \mathcal{X}$ will be
quasiconvex in $\mathcal Y$.

\begin{proof}[Proof of Theorem \ref{thm: engulf-alm-gferf}.] Note that $Q$ is a finite normal subgroup of $G$ by
Theorem \ref{thm:free-engulf}.

Consider the quotient $G_1=G/Q$ together with the natural homomorphism $\psi:G \to G_1$.
Since $Q$ is finite, $\psi$ is a quasiisometry between $G$ and $G_1$ ($G_1$ is equipped with the word metric
induced by the image of the finite generating set of $G$). Therefore, $G_1$ is also hyperbolic
and any preimage map ${\bar \psi}^{-1}: G_1 \to G$ (which maps an element of $G_1$ to some element of $G$
belonging to the corresponding left coset modulo $Q$) is a quasiisometry as well.

Choose an arbitrary proper quasiconvex subgroup $H_1\le G_1$. Then ${\bar \psi}^{-1}(H_1)$ is a quasiconvex subset
of $G$ and
$${\bar \psi}^{-1}(H_1) \subseteq \psi^{-1}(H_1) \subseteq {\bar \psi}^{-1}(H_1)\cdot Q,$$
where $\psi^{-1}(H_1)$ is the full preimage of $H_1$ in $G$.

As $Q$ is finite, the above formula implies
${\bar \psi}^{-1}(H_1) \approx \psi^{-1}(H_1)$. Therefore $\psi^{-1}(H_1)$ is quasiconvex in $G$ by
Remark \ref{rem:equiv-qc}. According to our assumptions, there is a proper finite index subgroup $L \le G$
containing $\psi^{-1}(H_1)$. By definition, $Q \le L$, hence $\psi(L)$ is a proper finite index subgroup of $G_1$
with $H_1 \le \psi(L)$.

Thus, we have shown that $G_1$ also engulfs each proper quasiconvex subgroup. By the construction, $G_1$ is
residually finite and, therefore, GFERF (Theorem \ref{thm:qc-englulf_2}).

Consider any quasiconvex subgroup $H \le G$. Then $\psi(H)$ is quasiconvex in $G_1$ and, thus,
it is closed in the profinite topology on $G_1$. The homomorphism $\psi$ is a continuous map if $G$ and $G_1$
are equipped with their profinite topologies, thus the full preimage $\psi^{-1}(\psi(H))=H\cdot Q$ is closed
in $G$. Obviously $H,Q \le K$ (where $K$ is the profinite closure of $H$ in $G$), hence $K=HQ$. Q.e.d.
\end{proof}

\section{Free products of GFERF groups}\label{sec:free-prod-gferf}
In the previous section we considered hyperbolic groups which engulf every proper quasiconvex subgroup.
Let us name them {\it QE-groups}, for brevity.

As it can be seen from Theorem \ref{thm: engulf-alm-gferf}, any QE-group $G$ is very close to being GFERF.
In fact, $G$ is quasiisometric to the quotient $G/E(G)$ which is GFERF by Corollary \ref{cor:engulf-rf} and
Theorem \ref{thm:qc-englulf_2}. Nevertheless, the answer to the question whether each QE-group is GFERF is
still unclear. Theorem \ref{thm: engulf-alm-gferf}
would yield a positive answer if a free product of any two QE-groups were a QE-group itself.
Unfortunately, the author is unable to prove this; actually, he
doubts if this is true in general.

However, the following statement, proved by R. Burns, can be used to show that a free product of GFERF-groups is,
again, GFERF:

\begin{lemma} {\rm (\cite[Thm. 1.1]{Burns})} \label{lem:burns} Suppose $G$ is a free product of its subgroups $G_i$
indexed by some set $I$, and let $H$ be a finitely generated subgroup. If for each $i \in I$, $g \in G$, the subgroup
$(H^g \cap G_i)$ is $G_i$-separable, then $H$ is $G$-separable.
\end{lemma}

Let $H$ be a  subgroup of a group $G$.
Suppose $\cA$  and $\cB$ are finite generating sets for $G$ and $H$ respectively and $|\cdot|_G$, $|\cdot|_H$ are
the corresponding length functions. Set $\hat c=\max\{|b|_G~:~b \in \cB\}$. Evidently, $|h|_G \le \hat c |h|_H$
for all $h \in H$.

$H$ is called {\it undistorted} in $G$ if
there exists a constant $c \ge 0$ such that $|h|_H \le c|h|_G$ for every $h \in H$.
In a hyperbolic group $G$, a finitely generated
subgroup is undistorted if and only if it is quasiconvex (\cite[Lemma 1.6]{hyp}).

\begin{proof}[Proof of Theorem \ref{thm:free-prod-gferf}.] It is well known that a free product of hyperbolic groups is
a hyperbolic group (see, for instance, \cite[1.34]{Ghys}). Thus, $G$ is hyperbolic.
Clearly, the subgroups $G_1$ and $G_2$ are undistorted in $G$; consequently, they are quasiconvex.

Choose an arbitrary quasiconvex subgroup $H \le G$, an element $g \in G$ and $i \in \{1,2\}$.
The subgroup $H^g$ is quasiconvex by Remark \ref{rem:qc-shifts}. Since the intersection of two
quasiconvex subgroups is quasiconvex (\cite[Prop. 3]{Short}), $(H^g \cap G_i)$ is quasiconvex in $G$.
Consequently, $(H^g \cap G_i)$ is undistorted in $G$, and, hence, it is undistorted in $G_i$. Thus
$(H^g \cap G_i)$ is $G_i$-separable because $G_i$ is GFERF.

According to Lemma \ref{lem:burns}, $H$ is $G$-separable. Q.e.d.
\end{proof}

\end{document}